\documentclass[a4paper, 12pt]{article}
\usepackage{adjustbox}
\usepackage{graphicx}
\usepackage{amsfonts}
\usepackage{amsbsy}
\usepackage{amssymb}
\usepackage[fleqn]{amsmath}
\usepackage{multicol}
\usepackage{longtable,ltxtable,booktabs}
\usepackage{blkarray}
\usepackage{afterpage}
\usepackage{float}
\usepackage{caption}
\usepackage{subcaption}
\usepackage{multirow}
\usepackage{pdflscape}
\usepackage{latexsym}
\usepackage{rotating}
\usepackage{enumerate}

\usepackage{mathtools}

\usepackage{epstopdf}
\usepackage{caption}
\usepackage{xcolor}
\usepackage[nodisplayskipstretch]{setspace}
\usepackage{amsmath}
\usepackage{longtable}
\usepackage[margin=.8 in]{geometry}

\usepackage{verbatim}

\usepackage{mathrsfs}
\usepackage{amsthm}
\usepackage{setspace}
 \usepackage{pdflscape}
\usepackage{colortbl}
\usepackage{latexsym}
\usepackage{amsfonts}
\usepackage{epsfig}
\usepackage{hyperref}
\usepackage{graphics}
\usepackage{float}
\theoremstyle{definition}
\newtheorem{thm}{Theorem}[section]

\newtheorem{defn}{Definition}[section]
\newtheorem{prop}{Proposition}[section]

\newtheorem{lemma}{Lemma}[section]

\newtheorem{cor}{Corollary}[section]

\theoremstyle{definition}

\title{ Some results on GVZ-groups with two character degrees}

\date{{}}\author{Nabajit Talukdar\thanks{Corresponding author.Email address : ntalukdar2000@yahoo.co.in}\\
	\small and\\
	Kukil Kalpa Rajkhowa\thanks{ E-mail address : kukilrajkhowa@yahoo.com}\\
	\small Department of Mathematics\\
	\small Cotton University\\
	\small Guwahati-781001, India}

\begin{document}
\maketitle

\begin{abstract}

We investigate the finite groups $G$ for which $\chi(1)^{2}=|G:Z(\chi)|$ for all characters $\chi \in Irr(G)$ and $|cd(G)|=2$, where $cd(G)=\{\chi(1)| \chi \in Irr(G)\}$. We call such a group a GVZ-group with two character degrees. We establish bijections between  the sets of characters of some  groups obtained from a GVZ-group with two character degrees. Additionally we obtain some alternate characterizations of a GVZ-group with two character degrees and we construct a GVZ-group having the character degree set $\{1,p\}$.
 \end{abstract}

  $\mathbf{2010\;Mathematics\;Subject\;Classification:}\;20C15$
	
	$\mathbf{Keywords\;and\;Phrases:}$ GVZ-groups, character degrees, generalized Camina pair

\section{Introduction}
In this paper, all groups  are finite. By $Irr(G)$, $Lin(G)$ and $nl(G)$  we denote the set of   complex irreducible characters, the set  of  complex linear irreducible characters and the set of  complex non-linear irreducible characters of the group $G$ respectively. By $c(G)$ we denote nilpotency class of $G$. For any two subgroups $H,K$ of a group $G$,  $[H,K]$ denotes the subgroup of $G$ generated by the elements of the form $[h,k]=h^{-1}k^{-1}hk$, where $h\in H, k\in K$ and $Z(G)$ denotes the centre of the group $G$. For any $\chi \in Irr(G)$, the kernel and  the centre of $\chi$ are defined by  $ker(\chi)=\{g\in G: \chi(g)=\chi(1)\}$ and $Z(\chi)=\{g\in G: |\chi(g)|=\chi(1)\}$ respectively. For the terminologies not defined here follow the Isaacs' book \cite{isaacs1994character} and Rotman's book \cite{rotman2012introduction}. Collecting the information about the structure of a group from the set of character degrees is an important direction of research in group theory. There are many research articles (\cite{heydari2011finite}, \cite{chillag2008finite}, \cite{fernandez2001groups},  \cite{lewis2020groups}, \cite{berkovich1992distinct}) attempting to answer this question by examining the groups whose irreducible characters have some distinctive properties like "the group has only two character degrees", "degrees of the non-linear irreducible characters are distinct", "the centres of the irreducible characters form a chain with respect to inclusion" etc. In \cite{nenciu2016nested}, A. Nenciu obtained various characterizations of a finite group $G$ for which $\chi(1)^{2}=|G:Z(\chi)|$ for all characters $\chi \in Irr(G)$ and $\{Z(\chi):\chi \in Irr(G)\}$
 is a chain . In \cite{lewis2020groups}, Mark L. Lewis obtained the characterizations of the groups for which $\{Z(\chi):\chi \in Irr(G)\}$
 is a chain with respect to inclusion and obtained some informations regarding the structures of these groups. In this paper we examine the GVZ-groups having two character degrees. The concept of GVZ-group first appeared in \cite{murai1994} under the name "groups of Ono type". The term  GVZ-group was introduced by A. Nenciu. 
 
 \begin{defn}
 \cite{nenciu2016nested}
A non-Abelian group $G$ is called a generalized VZ-group (GVZ for short) if for all $\chi\in Irr(G)$ we have $\chi(g)=0$ for all $g\in G\setminus Z(\chi)$.
 \end{defn}
Let $H$ be a subgroup of a group $G$. For any character $\lambda$ of $H$, $\lambda^{G}$ denotes the induced character on $G$. 

\begin{defn}
\cite{isaacs1994character}
Let $H\leq G$ and $\lambda$ be a character of $H$. Then $\lambda^{G}$, the induced character on $G$, is given by 
$\lambda^{G}(g)=\frac{1}{|H|}\sum_{x\in G} \lambda^{o}(xgx^{-1})$, where $\lambda^{o}$ is defined by $\lambda^{o}(h)=h$ if $h\in H$ and $\lambda^{o}(y)=0$ if $y \not \in H$.
\end{defn}

\begin{prop}
\cite{alperin1995group}
Let $H\leq G$ and $\lambda$ be a character of $H$. If $T$ is a transversal for $H$ in $G$, then for any $g\in G$ we have 
$\lambda^{G}(g)=\sum_{t\in T}\lambda^{o}(tgt^{-1})$.
\end{prop}

 By Corollary 2.30 \cite{isaacs1994character} we get that if $G$ is a GVZ-group then $\chi(1)^{2}=|G:Z(\chi)|$. If $G$ is a GVZ-group with two character degrees then it follows that $|Z(\chi)|$ is constant for all $\chi \in nl(G)$. In \cite{lewis2020groups} Mark L. Lewis obtained various results about the factor groups for a group where centres of the irreducible characters form a chain. In \cite{talukdar2024gvz} Talukdar et. al. initiatied the study of GVZ-group with two character degrees. In this paper, we  establish bijections between  the sets of characters of some  groups obtained from a GVZ-group with two character degrees. Additionally we obtain some alternate characterizations of a GVZ-group with two character degrees  and construct a GVZ-groups with two character degrees. We obtain the following results.

\begin{thm}
Let $G$ be a GVZ-group with two character degrees and   $\chi \in nl(G)$. Let $Irr_{*}(Z(\chi)/[Z(\chi),G])=\{\lambda \in Irr(Z(\chi)): [Z(\chi),G] \subseteq ker \lambda, G' \not \subseteq ker \lambda\}$.  Then
\begin{enumerate}[(i)]

\item  
$|\{\phi \in nl(G): Z(\phi)=Z(\chi)\}|
=|Z(\chi)|(\frac{1}{|[Z(\chi),G]|}-\frac{1}{|G'|})$.

\item If $\lambda \in Irr_{*}(Z(\chi)/[Z(\chi),G])$, then $\lambda ^{G}$ has a unique non-linear irreducible character of $G$ as a constituent;

 \item There is a bijection between the sets $Irr_{*}(Z(\chi)/[Z(\chi),G])$ and $nl(G/[Z(\chi),G])$;

\item There is a bijection between the sets $Irr_{*}(Z(\chi)/[Z(\chi),G])$ and $\{ \phi \in nl(G): Z(\phi)=Z(\chi)\}$.

\end{enumerate}
\end{thm}

 \begin{thm}
 Let $G$ be a group with two character degrees. Then the following are equivalent:
 \begin{enumerate}[(i)]

   \item $G$ is a GVZ-group.
   
   \item For any $\chi \in Irr(G)$, $x[Z(\chi),G]=Cl_{G}(x)$ if $x\in Z(\chi) \setminus \cup_{\phi \in nl(G), Z(\phi)\neq Z(\chi)}Z(\phi)$. 

   \end{enumerate}
 \end{thm}

\begin{prop}
  Let $p$ be an odd prime and  $G_{n}$ be the group defined by 
   $G_{n}=<\alpha, \alpha_{1}, \alpha_{2}, \ldots , \alpha_{n}, \beta_{1}, \beta_{2}, \ldots , \beta_{n}\ | \ [\alpha_{i}, \alpha]=\beta_{i}, \alpha^{p}=\alpha_{i}^{p}=\beta_{i}^{p}=1>$. Then $G_{n}$ is a GVZ-group of order $p^{2n+1}$ with $cd(G_{n})=\{1, p\}$. The centres of the non-linear irreducible characters of $G_{n}$ are $\{\alpha_{i_{1}},\alpha_{i_{2}},\ldots , \alpha_{i_{n-1}}, \beta_{1}, \beta_{2}, \ldots , \beta_{n}\}$ where $\{ \alpha_{i_{1}}, \alpha_{i_{2}}, \ldots , \alpha_{i_{n-1}}\}$ varies over all the subsets of order $n-1$ of $\{\alpha_{1}, \alpha_{2}, \ldots , \alpha_{n}\}$.

\end{prop}

\section{Preliminaries}

The aim of this section is to prove some essential results which will be applied in the proofs of our main results. We state the following results.

\begin{lemma}
\label{lemma1}
\cite{talukdar2024gvz}
Let $\chi \in Irr(G)$. Then $\chi$ is a linear character if and only if $[Z(\chi), G]=G'$.
\end{lemma}

\begin{lemma}
\label{lemma2}
\cite{talukdar2024gvz}
Let $\chi \in Irr(G)$. Then $Z(G/[Z(\chi),G])=Z(\chi)/[Z(\chi), G]$.
\end{lemma}

The following results can be found in \cite{isaacs1994character}.

\begin{lemma}
\label{lemma3}
\cite{isaacs1994character}
Let $H\leq G$ and $\chi$ be a character of $G$. Then $[\chi_{H}, \chi_{H}]\leq [G:H][\chi,\chi]$ with equality if and only if $\chi(g)=0$ for all $g\in G\setminus H$.
\end{lemma}



     

\begin{lemma}
\label{lemma4}
\cite{isaacs1994character}
Let $\chi \in Irr(G)$. Then $\chi(1)^{2}\leq |G/Z(\chi)|$. Equality holds if and only if $\chi$ vanishes on $G \setminus Z(\chi)$.
\end{lemma}

\begin{lemma}
\label{lemma5}
\cite{isaacs1994character}
Let $\chi \in Irr(G)$ and $|G/Z(\chi)|$ is Abelian. Then  $\chi(1)^{2}=|G/Z(\chi)|$.
\end{lemma}

\begin{lemma}
\label{lemma6}
\cite{isaacs1994character}
Let $\chi$ be an irreducible character of a group $G$. Then $Z(\chi)/ker \chi=Z(G/ker(\chi))$.
\end{lemma}

\begin{lemma}
\label{lemma7}
\cite{isaacs1994character}
Let $\chi$ be an irreducible character of a group $G$ and $N\unlhd G$. Then $\chi \in Irr(G/N)$ if and only if $N \subseteq ker \chi$.
\end{lemma}

From Lemma \ref{lemma6} we get the following result.

\begin{cor}
\label{cor1}
    Let $G$ be a group and $\chi \in \text{Irr}(G)$. Then $[Z(\chi), G]\leq ker \chi$.
\end{cor}

In the following lemma we prove a necessary and sufficient condition when the centres of two non-linear irreducible characters of a GVZ-group with two character degrees are equal.

\begin{lemma}
\label{lemma8}
Let $G$ be a GVZ group with two character degrees. For $\chi, \phi \in nl(G)$ , $Z(\chi)=Z(\phi)$ if and only if $\phi \in nl(G/[Z(\chi), G])$.

\begin{proof}
First we assume that $Z(\chi)=Z(\phi)$. Then $[Z(\chi),G]=[Z(\phi),G]\leq ker \phi$ and hence $\phi \in nl(G/[Z(\chi), G])$. Conversely if $\phi \in nl(G/[Z(\chi), G])$ then $[Z(\chi), G] \leq ker \phi$ and hence $Z(\chi)=Z(\phi)$.

\end{proof}

\end{lemma}

From Lemma \ref{lemma2} and Lemma \ref{lemma7}, we can deduce when a centre of a non-linear irreducible character is a subset of a centre of a another irreducible character.

\begin{cor}
\label{cor2}
 Let $G$ be a group and $\chi$ be a non-linear irreducible character of $G$. Then for any $\phi \in Irr(G)$, $Z(\chi)\subseteq Z(\phi)$ if and only if $\phi \in Irr(G/[Z(\chi),G])$.
 \begin{proof}
First suppose that  $Z(\chi)\subseteq Z(\phi)$. Then $[Z(\chi),G]\subseteq [Z(\phi),G] \subseteq ker\phi$ and hence $\phi \in Irr(G/[Z(\chi),G])$. Conversely assume that $\phi \in Irr(G/[Z(\chi),G])$. Since $Z(G/[Z(\chi),G])=Z(\chi)/[Z(\chi), G]$ we get that 
$Z(\chi)\subseteq Z(\phi)$.
 \end{proof}
\end{cor}

A pair $(G,N)$ is said to be a generalized Camina pair (abbreviated GCP) if $N$ is normal in $G$ and all the non-linear irreducible characters of $G$ vanish outside $N$. 
The notion of GCP was introduced by Lewis in \cite{lewis2009vanishing}. An equivalent condition for a pair $(G,N)$ to be a GCP is: A pair $(G,N)$ is a GCP if and only if for $g\in G\setminus N$, the conjugacy class of $g$ in $G$ is $gG'$.
In \cite{prajapati2017irreducible} the authors obtained the following results.

\begin{thm}
\label{thm1}

Let $(G,Z(G))$ be a GCP. Then we have the following.

   \begin{enumerate}[(1)]

       \item $cd(G)=\{1, |G/Z(G)|^{\frac{1}{2}}\}$.

       \item The number of non-linear irreducible characters of $G$ is $|Z(G)|-|Z(G)/G'|$.

   \end{enumerate}
\end{thm}





Now we get a count of the number of number of non-linear irreducible characters of certain types of GVZ-groups with two character degrees.
\begin{lemma}
\label{lemma9}
Let $G$ be a GVZ-group with two character degrees. Suppose the centres of all the non-linear irreducible characters are equal. Then the number of non-linear characters of $G$ is $|Z(G)|-|Z(G)/G'|$.

\begin{proof}

We get that $Z(\chi)=Z(G)$ and $\chi^{2}(1)=|G/Z(\chi)|=|G/Z(G)|$ for all $\chi \in nl(G)$. Thus $(G,Z(G))$ is a generalized camina pair and hence the number of non-linear irreducible characters is $|Z(G)|-|Z(G)/G'|$.
\end{proof}
\end{lemma}

\begin{lemma}
\label{lemma10}
Let $G$ be a GVZ-group with two character degrees and   $\chi \in nl(G)$. Then
\begin{enumerate}[(i)]

\item $cd(G)=\{1, |G/Z(\chi)|^{\frac{1}{2}}\}$.

\item The number of non-linear irreducible characters of $G$ is
$|Z(\chi)|-|Z(\chi)/G'|$.


\end{enumerate}
\begin{proof}

\begin{enumerate}[(i)]

\item Obvious.

\item For all $\chi \in nl(G)$, $\chi(1)^{2}=|G/Z(\chi)|$. Let $t$ be the number of non-linear irreducible characters of $G$. Then $|G|=|G/G'|+t|G/Z(\chi)|$. Hence $t=|Z(\chi)|-|Z(\chi)/G'|$.


\end{enumerate}
\end{proof}
\end{lemma}

In order to construct a GVZ-group with the set of character degrees $\{1,p\}$, we rely on the classification of groups of orders $p^{4}$, $p^{5}$ and $p^{6}$ done by Rodney James in \cite{james1980groups}. In \cite{prajapati2017irreducible} Prajapati et al. found many character related results of the groups listed in \cite{james1980groups}. Here we cite the following result from \cite{prajapati2017irreducible}.

\begin{prop}
\label{prop1}
\cite{prajapati2017irreducible}
Let $G$ be a non-Abelian $p$-group of order $p^{4}$. Then $(G,Z(G))$ is GCP if and only if nilpotence class of $G$ is $2$.

\end{prop}

\section{Main Results}

   

   


\begin{thm}
Let $G$ be a GVZ-group with two character degrees and   $\chi \in nl(G)$. Let $Irr_{*}(Z(\chi)/[Z(\chi),G])=\{\lambda \in Irr(Z(\chi)): [Z(\chi),G] \subseteq ker \lambda, G' \not \subseteq ker \lambda\}$.  Then
\begin{enumerate}[(i)]

\item  
$|\{\phi \in nl(G): Z(\phi)=Z(\chi)\}|
=|Z(\chi)|(\frac{1}{|[Z(\chi),G]|}-\frac{1}{|G'|})$.

\item If $\lambda \in Irr_{*}(Z(\chi)/[Z(\chi),G])$, then $\lambda ^{G}$ has a unique non-linear irreducible character of $G$ as a constituent;

 \item There is a bijection between the sets $Irr_{*}(Z(\chi)/[Z(\chi),G])$ and $nl(G/[Z(\chi),G])$;

\item There is a bijection between the sets $Irr_{*}(Z(\chi)/[Z(\chi),G])$ and $\{ \phi \in nl(G): Z(\phi)=Z(\chi)\}$.

\end{enumerate}

 \begin{proof}

 \begin{enumerate}[(i)]
\item
Let $\chi \in nl(G)$ and $\overline{G}=G/[Z(\chi),G]$. By Lemma  \ref{lemma8} we get that if $Z(\phi)=Z(\chi)$ if and only if $\phi \in nl(\overline{G})$. The group $\overline{G}$ is a GVZ group with two character degrees where the centres of all the non-linear irreducible characters are equal. 
Hence 
$ |\{\phi \in nl(G): Z(\phi)=Z(\chi)\}|
= |nl(\overline{G})|
=|Z(\overline{G})|-|Z(\overline{G})/\overline{G}'|
 =|Z(\chi)|(\frac{1}{|[Z(\chi),G]|}-\frac{1}{|G'|})
$

\item  Let $\lambda \in Irr_{*}(Z(\chi)/[Z(\chi),G])$. If all of the irreducible constituents of $\lambda ^{G}$ is linear, then we get that $G'\subseteq ker\ \lambda ^{G} \subseteq ker\lambda$, a contradiction. Let $\theta$ be a non-linear irreducible constituent of $\lambda  ^{G}$. Now 
$[Z(\chi),G] \subseteq \cap_{g\in G} (ker \lambda)^{g} \subseteq ker \theta$. Thus $\theta \in Irr(G/[Z(\chi),G])$. Hence by Corollary \ref{cor2} we get that $Z(\theta)=Z(\chi)$. This gives that $\theta$ vanishes outside $Z(\chi)$. From the fact that $\theta_{Z(\theta)}=e\lambda$ and $\theta(1)=|G/Z(\chi)|^{\frac{1}{2}}$, it follows that  $e=\frac{|G/Z(\chi)|^{\frac{1}{2}}}{\lambda(1)}$. Thus we get that 

\[
\theta(g)=
\begin{cases}
0 & \text{if} \ g \not \in Z(\chi),\\
\frac{|G/Z(\chi)|^{\frac{1}{2}}}{\lambda(1)}\lambda & \text{if} \ g \in Z(\chi).\\
  
\end{cases}
\]
This proves that $\theta$ is the unique non-linear irreducible constituent of $\lambda  ^{G}$.

\item Let $\lambda \in Irr_{*}(Z(\chi)/[Z(\chi),G])$. By the proof of (ii) above we get that $\lambda ^{G}$ has a unique non-linear irreducible constituent $\theta \in Irr (G/[Z(\chi), G])$. The mapping $\lambda \rightarrow \theta$ defines a one-one mapping between the sets $Irr_{*}(Z(\chi)/[Z(\chi),G])$ and $nl(G/[Z(\chi),G])$. Again let $\theta \in Irr(G/[Z(\chi), G])$. By Corollary \ref{cor2} it follows that $Z(\theta)=Z(\chi)$. Thus $\theta$ vanishes outside $Z(\chi)$. Hence we get that 
\[
\theta(g)=
\begin{cases}
0 & \text{if} \ g \not \in Z(\chi),\\
e\lambda & \text{if} \ g \in Z(\chi).\\
  
\end{cases}
\] This shows that the mapping $\lambda \rightarrow \theta$ is onto.

\item This follows from (iii) above since $Z(\chi)=Z(\theta)$ if and only if $\theta \in nl(G/[Z(\chi), G])$.
 
 \end{enumerate}

 \end{proof}
\end{thm}

In \cite{nenciu2016nested} A. Nenciu obtained some characterizations of nested GVZ-group. In the following theorem we obtain similar results in case of GVZ-group with two character degrees.

\begin{thm}
 Let $G$ be a non Abelian group with two character degrees. Then the following statements are equivalent:
 \begin{enumerate}[(i)]

   \item $G$ is a GVZ-group.
   
   \item For any $\chi \in Irr(G)$, $x[Z(\chi),G]=Cl_{G}(x)$ if $x\in Z(\chi) \setminus \cup_{\phi \in nl(G), Z(\phi)\neq Z(\chi)}Z(\phi)$. 


 \end{enumerate}

 \begin{proof}
    $(i) \Rightarrow (ii)$
        Assume that $G$ is a GVZ-group with two character degrees. First we consider a linear character $\chi \in Irr(G)$ so that $Z(\chi)=G$. Suppose $x\in G\setminus \cup_{\phi \in nl(G)}Z(\phi)$.
        Then $|C_{G}(x)|=\sum_{\chi \in Lin(G)}|\chi(x)|^{2}+ \sum_{\phi \in nl(G)}|\phi(x)|^{2}$
        $=|G/G'|+0=|G/G'|$.\\
        Hence $|Cl_{G}(x)|=|G/C_{G}(x)|=|G'|=|xG'|$. Since $Cl_{G}(x) \subseteq xG'$, we get that $x[Z(\chi),G]=xG'=Cl_{G}(x)$.     
        Next we consider $\chi \in nl(G)$ and $x\in Z(\chi) \setminus \cup_{\phi \in nl(G), Z(\phi)\neq Z(\chi)}Z(\phi)$. Then $Cl_{G}(x)\subseteq x[Z(\chi),G]$. We have
       \begin{align*}
          |C_{G}(x)| & = \sum_{\chi \in Irr(G)}|\chi(x)|^{2}&\\
          &=\sum_{\chi \in Irr(G)}\chi(1)^{2}&\\
          &= |G/G'|+ \sum_{\phi \in nl(G), Z(\phi) = Z(\chi)}|G/Z(\chi)|&\\
          &=|G/G'|+ |Z(\chi)|(\frac{1}{|[Z(\chi),G]|}-\frac{1}{|G'|})|G/Z(\chi)| &\\
          &=|G/[Z(\chi),G]|&
       \end{align*}
 We know that $|C_{G}(x)|=|\frac{G}{Cl_{G}(x)}|$. Hence it follows that $Cl_{G}(x)=x[Z(\chi),G]$.\\
$(ii) \Rightarrow (i)$

Let $\chi\in nl(G)$. Suppose $Z(\chi)\neq Z(G)$. Then we consider the group $\overline{G}=G/[Z(\chi),G]$. We get by induction on $|G|$ that $\overline{G}$ is a GVZ-group. Hence $\chi(1)^{2}=|G/Z(\chi)|$. Also by Corollary \ref{cor2} we get that if $\chi \in nl(G)$ and  $Z(\chi)\neq Z(G)$ then 
$|\{\phi \in nl(G): Z(\phi)=Z(\chi)\}| =|\{\phi \in nl(\overline{G}): Z(\phi)=Z(\chi)\}|=
|Z(\chi)|(\frac{1}{|[Z(\chi),G]|}-\frac{1}{|G'|})$.\\

Next we consider the case that there is a $\chi \in nl(G)$ such that $Z(\chi)= Z(G)$. We choose $x\in G \setminus Z(G)$.The following two cases arise :\\
\\
Case I: There is a $\theta \in nl(G)$ such that $Z(\theta)\neq Z(G)$ and $x\in Z(\theta)\setminus Z(G)$.\\
The number of $\theta \in nl(G)$ such that $x\in Z(\theta)\setminus Z(G)$ is at least $|Z(\theta)|(\frac{1}{|[Z(\theta),G]|}-\frac{1}{|G'|})$. Moreover, $\theta(1)^{2}=|G/Z(\theta)|$. Suppose $\psi \in nl(G)$ be such that $x\in Z(\psi) \setminus Z(G)$ but $Z(\theta) \not \subseteq Z(\psi)$. Since $G$ has two character degrees we get that $\psi(1)^{2}=|G/Z(\psi)|=\theta(1)^{2}=|G/Z(\theta)|$ and hence $|Z(\psi)|=|Z(\theta)|$. Then we get that $x[Z(\theta),G]=Cl_{G}(x)$ .\\
Now $|G/[Z(\theta),G]|$\\
$=|G/Cl_{G}(x)|$\\
$=|C_{G}(x)|$\\
$=\sum_{\phi \in Irr(G)}|\phi(x)|^{2}$\\
$=|G/G'|$+ $|G/Z(\theta)||Z(\theta)|(\frac{1}{|[Z(\theta),G]|}-\frac{1}{|G'|})+$ $ \sum_{x\in Z(\psi) \setminus Z(G)}|G/Z(\psi)||Z(\psi)|(\frac{1}{|[Z(\psi),G]|}-\frac{1}{|G'|})$ + $\sum_{\chi \in nl(G), Z(\chi)=Z(G)}|\chi(x)|^{2}$\\
$=|G/[Z(\theta),G]|+  \sum_{x\in Z(\psi) \setminus Z(G)}|G/Z(\psi)||Z(\psi)|(\frac{1}{|[Z(\psi),G]|}-\frac{1}{|G'|}) + \sum_{\chi \in nl(G), Z(\chi)=Z(G)}|\chi(x)|^{2}$\\
$\Rightarrow \sum_{x\in Z(\psi) \setminus Z(G)}|G/Z(\psi)||Z(\psi)|(\frac{1}{|[Z(\psi),G]|}-\frac{1}{|G'|}) + \sum_{\chi \in nl(G), Z(\chi)=Z(G)}|\chi(x)|^{2}=0$.\\
This gives that $\chi(x)=0$ for all $\chi \in nl(G)$ such that $Z(\chi)= Z(G)$.\\
\\
Case II: We consider $x\not \in Z(\theta)$ for all $\theta \in nl(G)$.\\
 In this case we get that $xG'=Cl_{G}(x)$.\\ 
 Thus $|G/G'|=|G/Cl_{G}(x)|=|C_{G}(x)|$\\
 $=|G/G'|+\sum_{\theta \in nl(G), Z(\theta) \neq Z(G)}|\theta(x)|^{2}+ \sum_{\chi \in nl(G), Z(\chi) = Z(G)}|\chi(x)|^{2}$\\
 $=|G/G'|+0+ \sum_{\chi \in nl(G), Z(\chi) = Z(G)}|\chi(x)|^{2}$.\\
 This gives that $\chi(x)=0$ for all $\chi \in nl(G)$ such that $Z(\chi)= Z(G)$.\\

    

 \end{proof}
\end{thm}

In the following Proposition we construct a GVZ-group having the character degree set $\{1,p\}$.

\begin{prop}
   Let $p$ be an odd prime and  $G_{n}$ be the group defined by 
   $G_{n}=<\alpha, \alpha_{1}, \alpha_{2}, \ldots , \alpha_{n}, \beta_{1}, \beta_{2}, \ldots , \beta_{n}\ | \ [\alpha_{i}, \alpha]=\beta_{i}, \alpha^{p}=\alpha_{i}^{p}=\beta_{i}^{p}=1>$. Then $G_{n}$ is a GVZ-group of order $p^{2n+1}$ with $cd(G_{n})=\{1, p\}$. The centres of the non-linear irreducible characters of $G_{n}$ are $\{\alpha_{i_{1}},\alpha_{i_{2}},\ldots , \alpha_{i_{n-1}}, \beta_{1}, \beta_{2}, \ldots , \beta_{n}\}$ where $\{ \alpha_{i_{1}}, \alpha_{i_{2}}, \ldots , \alpha_{i_{n-1}}\}$ varies over all the subsets of order $n-1$ of $\{\alpha_{1}, \alpha_{2}, \ldots , \alpha_{n}\}$.
   \begin{proof}
       We get that $G_{n}'=Z(G_{n})=<\beta_{1}, \beta_{2}, \ldots , \beta_{n}>$. 
       We prove by induction. For $n=1$, $G_{1}=<\alpha, \alpha_{1}, \beta_{1}| [\alpha_{1}, \alpha]=\beta_{1}, \alpha^{p}=\alpha_{1}^{p}=\beta_{1}^{p}=1>$
       is a group of order $p^{3}$. Since all the non-linear irreducible characters of $G_{1}$ are faithful (see \cite{gorenstein2007finite}), it follows that the centres of the non-linear irreducible characters of $G_{1}$ are $Z(G_{1})=\{\beta_{1}\}$. For $n=2$,  $G_{2}=<\alpha, \alpha_{1}, \alpha_{2}, \beta_{1}, \beta_{2}| [\alpha_{i}, \alpha]=\beta_{i}, \alpha^{p}=\alpha_{i}^{p}=\beta_{i}^{p}=1 (i=1, 2)> $ is the group $\phi_{4}(1^{5})$ in \cite{james1980groups}. Since $|G_{2}/Z(G_{2})|=p^{3}$, we get that $cd(G_{2})=\{1,p\}$. The proper normal subgroups of order $p^{3}$ of $G_{2}$ containing $Z(G_{2})$ are $N_{1}=\{\alpha_{1}, \beta_{1}, \beta_{2}\}$, 
        $N_{2}=\{\alpha_{2}, \beta_{1}, \beta_{2}\}$ and $N_{3}=\{\alpha, \beta_{1}, \beta_{2}\}$.
       Consider the normal subgroups $K_{1}=<\beta_{1}>$ and $K_{2}=<\beta_{2}>$ of $G_{2}$. Then
        $G_{2}/K_{1}=<\overline{\alpha}, \overline{\alpha_{1}}, \overline{\alpha_{2}}, \overline{\beta_{2}}|
[\overline{\alpha_{2}}, \alpha]=\overline{\beta_{2}}, \overline{\alpha}^{p}=\overline{\alpha}_{1}^{p}=\overline{\alpha}_{2}^{p}=\overline{\beta}_{2}^{p}=1>$  and $Z(G_{2}/K_{1})=\{ \overline{\alpha_{1}}, \overline{\beta_{2}}\}$. By Proposition \ref{prop1}, $(G_{2}/K_{1}, Z(G_{2}/K_{1}))$ is a GCP. Hence for any $\chi \in nl(G_{2}/K_{1})$, $Z(\chi)=\{\alpha_{1}, \beta_{1}, \beta_{2}\}$. Similarly for any $\chi \in nl(G_{2}/K_{2})$, $Z(\chi)=\{\alpha_{2}, \beta_{1}, \beta_{2}\}$. Since $[N_{3}, G_{2}]=G_{2}'$,  by Lemma \ref{lemma1}, $N_{3}$ can not be the centre of any non-linear irreducible character. By induction hypothesis we assume that $G_{n-1}$  is a GVZ-group of order $p^{2n-1}$ with $cd(G_{n-1})=\{1, p\}$ and the centres of the non-linear irreducible characters of $G_{n-1}$ are $\{\alpha_{i_{1}},\alpha_{i_{2}},\ldots , \alpha_{i_{n-2}}, \beta_{1}, \beta_{2}, \ldots , \beta_{n-1}\}$ where  $\{ \alpha_{i_{1}}, \alpha_{i_{2}}, \ldots , \alpha_{i_{n-2}}\}$ vary over all subsets of order $n-2$ of $\{\alpha_{1}, \alpha_{2}, \ldots , \alpha_{n-1}\}$. Consider $\chi \in nl(G_{n})$ and $\beta_{j}\in Z(G_{n}) \cap ker(\chi)$. For $K_{j}=<\beta_{j}>\unlhd G_{n}$,
$G_{n}/K_{j}=<\overline{\alpha}, \overline{\alpha_{1}}, \overline{\alpha_{2}},
\ldots , \overline{\alpha_{n}},\overline{\beta_{1}},\overline{\beta_{2}}, \ldots , \overline{\beta_{j-1}},\overline{\beta_{j+1}},\ldots , \overline{\beta_{n}}| [\overline{\alpha_{i}}, \overline{\alpha}]=\overline{\beta_{i}}, \overline{\alpha}^{p}=\overline{\alpha_{i}}^{p}=\overline{\beta_{i}}^{p}=1>$
$\cong G_{n-1} \times <\overline{\alpha_{j}}>$. We get that $\chi \in nl(G_{n}/K_{j})$ and $Z(\chi)=Z(\theta)\times <\overline{\alpha_{j}}>$ for some $\theta \in nl(G_{n-1})$. Thus $Z(\chi)= <\overline{\alpha_{i_{1}}}, \overline{\alpha_{i_{2}}},
\ldots , \overline{\alpha_{i_{n-2}}},\overline{\beta_{1}},\overline{\beta_{2}}, \ldots , \overline{\beta_{j-1}},\overline{\beta_{j+1}},\ldots , \overline{\beta_{n}}> \times <\overline{\alpha_{j}}>$, where $\{ \overline{\alpha_{i_{1}}}, \overline{\alpha_{i_{2}}},
\ldots , \overline{\alpha_{i_{n-2}}}\}$ is a subset of $\{ \overline{\alpha_{1}}, \overline{\alpha_{2}},
\ldots , \overline{\alpha_{j-1}}, \overline{\alpha_{j+1}}, \ldots , \overline{\alpha_{n}} \}$. Thus for $\chi \in nl(G_{n})$,
$Z(\chi)=<\alpha_{i_{1}},\alpha_{i_{2}},\ldots , \alpha_{i_{n-2}}, \alpha_{j}, \beta_{1}, \beta_{2}, \ldots , \beta_{n}>$.
       
   \end{proof}
\end{prop}

\end{document}